\documentclass[a4paper,11pt]{article}
\usepackage{jheppub} % for details on the use of the package, please see the JHEP-author-manual
\pdfoutput=1
\usepackage[utf8]{inputenc}
\usepackage{amsmath,amsthm,amssymb,amsfonts,mathtools,bm,array,float}
\usepackage{caption,subcaption}
\usepackage{graphicx}
\usepackage[T1]{fontenc} 

\numberwithin{equation}{section}
\numberwithin{figure}{section}
\numberwithin{table}{section}

\begin{document}

\title{PINNs in More General Geometry}
\author[a]{Edward Hirst}
\affiliation[a]{Instituto de Matemática, Estatística e Computação Científica (IMECC), Universidade Estadual de Campinas (UNICAMP), 13083-859, Brazil}
\emailAdd{ehirst@unicamp.br}

%\preprint{\begin{flushright}
%...
%\end{flushright}}
\abstract{
Neural architectures trained with losses inspired by differential conditions are the basis for PINN models. Since many constructions in differential geometry may be framed as minimisation of a differential functional, these functionals can be coded as loss functions to align the AI loss-minimisation goal with that of solving the geometric problem.
This contribution to the Recent Progress in Computational String Geometry workshop proceedings introduces the PINN architecture defining principles, motivates how they are well suited for problems in differential geometry, and demonstrates their use via summaries of three works at this intersection.
}

\maketitle

%%%%%%%%%%%%%%%%%%%%%%%%%%%%%%%%%%%%%%%%%%%%%%%%%%%%
\section{Introduction}\label{sec:intro}
Many central problems in differential geometry ask for distinguished metrics, curvatures, or immersions selected by differential constraints or variational principles. Einstein geometry seeks metrics satisfying $Ric(g)=\lambda g$; the Nirenberg problem asks which functions on $S^2$ occur as Gaussian curvature in a fixed conformal class; and Willmore theory studies critical points of the bending functional $\mathcal{W}(\Sigma)=\int_\Sigma H^2\,dA$ \cite{Besse:1987pua,Nirenberg1953,Willmore1965}. In each case the unknown object is geometric and global, but the condition defining it is local and differential. This makes geometry a natural arena for numerical functional minimisation, one chooses an ansatz for the unknown field, samples points on the underlying manifold, and drives a residual or energy functional toward zero while respecting global compatibility conditions.

Modern machine learning follows the same broad pattern. A neural network is a parametrised function, and training amounts to optimising parameters so as to minimise a differentiable objective. In the physics-informed neural network (PINN) paradigm, the loss is not built from labelled output data alone; it is designed directly from the governing equations, variational principles, and consistency conditions of the problem \cite{Raissi2019,Karniadakis2021}. Efficient automatic differentiation supplies exact derivatives of the network output with respect to both inputs and trainable parameters, so geometric quantities such as Christoffel symbols, Ricci tensors, Laplacians, or mean curvature can be evaluated without introducing a mesh-based finite-difference stencil \cite{Baydin2018}.

This structural alignment makes differential geometry particularly well suited to PINN methods. Manifolds provide an effectively inexhaustible supply of sample points, geometric constraints are often naturally expressible as differentiable residuals, and the learnt model is itself a smooth function that interpolates between sampled locations. Recent work has already shown that this viewpoint is productive across several geometric settings, including Einstein metrics on spheres \cite{Hirst:2025seh, Cortes:2026kfx}, Calabi--Yau manifolds \cite{douglas2020numerical, Larfors:2021pbb, Berglund:2022gvm}, and $\mathrm{G}_2$-manifolds \cite{douglas2024harmonic, Heyes:2026rch}, some minimal surface learning \cite{ZhouYe2023MinimalSurfacePINN, Hashimoto:2025zmi}; with examples where trained AI models have then led to analytically rigorous computer-assisted proofs \cite{GomezSerranoSurvey, wang2023, platt_nirenberg}. 
The goal of this proceedings article is to summarise the core principles behind these architectures and then illustrate them through three recent case studies in increasingly general geometric settings \cite{Hirst:2025seh,Cortes:2026kfx,Hirst:2026qwi}.

%%%%%%%%%%%%%%%%%%%%%%%%%%%%%%%%%%%%%%%%%%%%%%%%%%%%
\section{PINNs}\label{sec:pinns}
A feedforward neural network is a parametrised map built from affine layers and nonlinear activations,
\begin{equation}
f_\theta(x)=W_L\,\sigma\Bigl(W_{L-1}\,\sigma\bigl(\cdots \sigma(W_1x+b_1)\cdots\bigr)+b_{L-1}\Bigr)+b_L,
\label{eq:nn_general}
\end{equation}
where $\theta=\{W_\ell,b_\ell\}$ denotes the trainable parameters. The universal approximation theorems show that even relatively simple networks can approximate continuous functions on compact domains to arbitrary accuracy \cite{Cybenko1989,Hornik1991}. For geometry this matters because the unknown being sought is usually already a smooth function or tensor field; in these examples: a metric in local coordinates, a conformal factor on $S^2$, or an immersion into Euclidean space. The architectural question is therefore not merely expressive power, but how best to encode the geometric object so that some of its defining properties are built in from the outset.

In practice this means matching the network input and output spaces to the geometry. If the manifold is naturally described by multiple coordinate patches, one may use an \emph{Atlas Architecture}: a separate subnetwork learns the object on each patch, and additional losses enforce consistency on overlaps. If instead the geometry admits a convenient global embedding, one may use an \emph{Embedding Architecture}: the network learns a single global function on that embedded domain and global compatibility becomes automatic. 
Spectral feature maps, convolutional layers, attention blocks, or hybrid architectures may then be inserted when they better respect the symmetries or scales of the target problem. The guiding principle is always the same, reserve network capacity for the genuinely unknown part of the geometry and hard-code as much structure into the architecture to satisfy necessary geometric conditions intrinsically.

\begin{figure}[b]
\centering
\begin{subfigure}[t]{0.48\textwidth}
\centering
\includegraphics[width=0.95\textwidth]{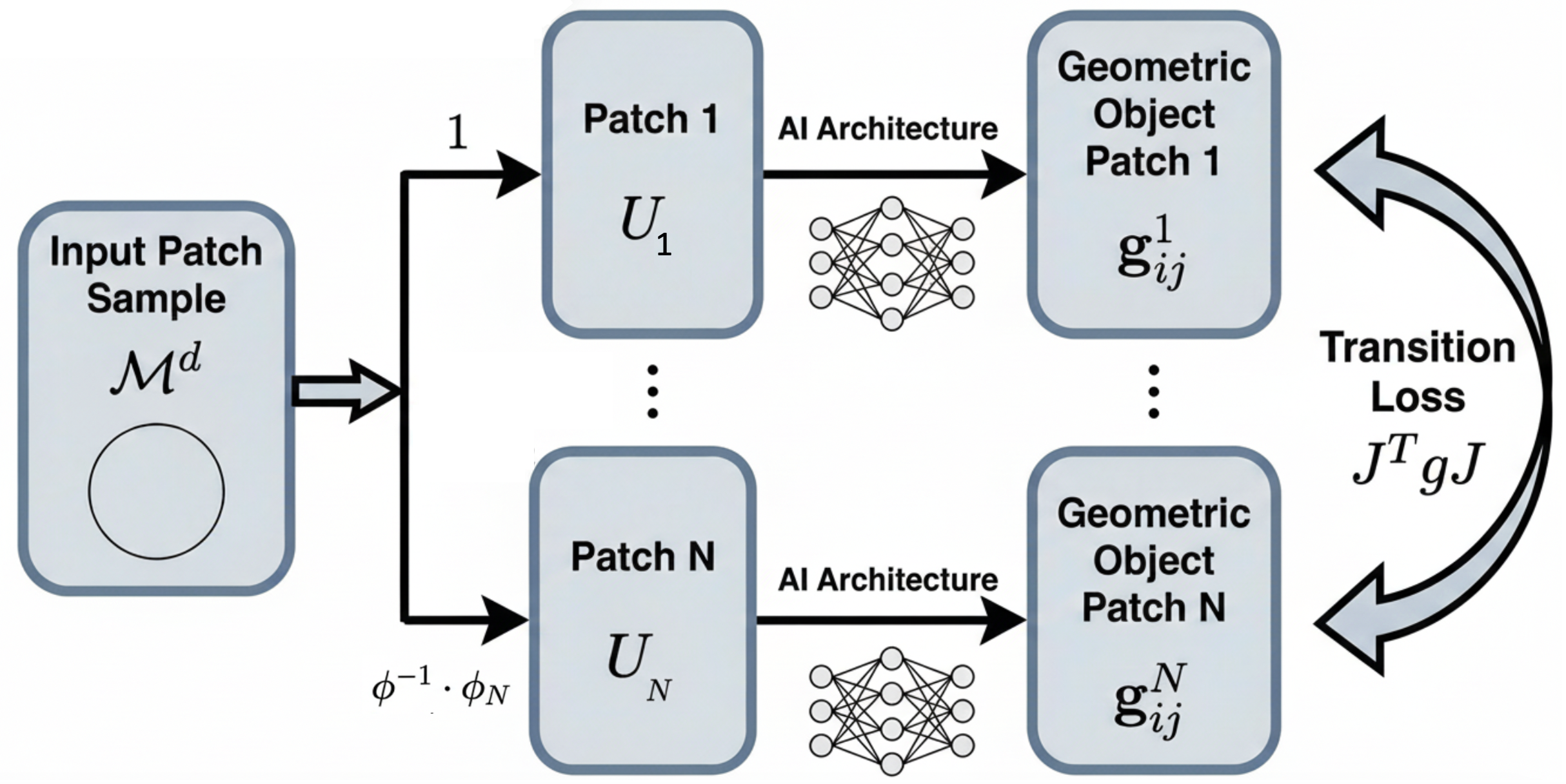}
\caption{Atlas Architecture: local subnetworks plus overlap consistency.}
\label{fig:atlas_arch}
\end{subfigure}
\hfill
\begin{subfigure}[t]{0.48\textwidth}
\centering
\includegraphics[width=0.95\textwidth]{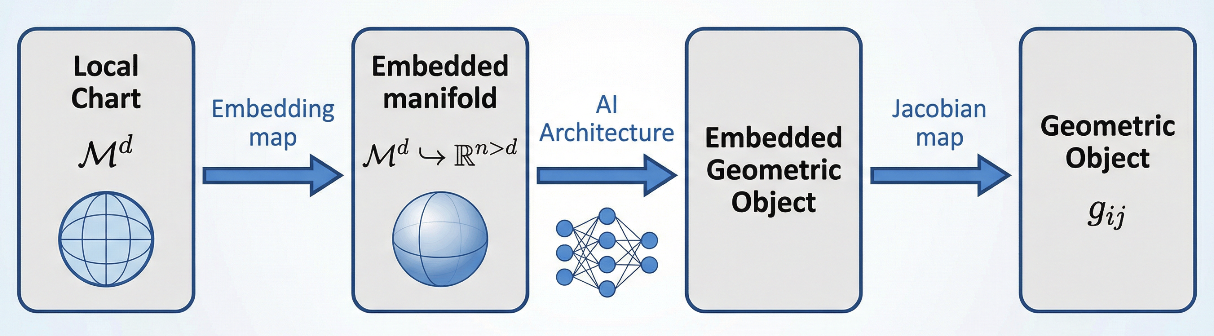}
\caption{Embedding Architecture: one global model on a known embedding.}
\label{fig:embedded_arch}
\end{subfigure}
\caption{Two recurrent design patterns for geometric PINNs. Atlas models trade a more local representation for explicit transition losses, while embedding models avoid overlap terms when a suitable global realisation is available.}
\label{fig:pinn_architectures}
\end{figure}

What makes a model a PINN is the loss. Instead of supervising against a labelled target field, one optimises a weighted combination of physics- or geometry-derived terms.
A PINN loss may be considered to take the following general form:
\begin{equation}
\mathcal{L}_{\mathrm{PINN}} = \lambda_{\mathrm{diff}}\, \mathcal{L}_{\mathrm{diff}} + \lambda_{\mathrm{bc}}\, \mathcal{L}_{\mathrm{bc}} + \lambda_{\mathrm{rep}}\, \mathcal{L}_{\mathrm{rep}}\,,
\label{eq:generic_pinn_loss}
\end{equation}
where $\mathcal{L}_{\mathrm{diff}}$ denotes \emph{differential losses}, and the central optimisation goal, built from PDE or Euler--Lagrange residuals.
In these examples these are constructed to reflect  the geometric defining properties from $Ric(g)-\lambda g$, $1-\Delta_{g_0}u-Ke^{2u}$, or the mean-curvature integrand entering Willmore energy. 
The term $\mathcal{L}_{\mathrm{bc}}$ stands for \emph{boundary, overlap, or constraint losses}: patch transition conditions, periodic identifications, or smooth gluing relations all fall into this class. 
Finally $\mathcal{L}_{\mathrm{rep}}$ denotes \emph{repeller or regularity losses}, used to discourage collapse to degenerate metrics, singular embeddings, already-known undesired solutions, or numerically unstable regions of parameter space. 
In many geometric problems these regularisers are essential, because a differential residual alone does not prevent the optimiser from drifting toward pathological but low-loss configurations.

The quality of a geometric PINN depends at least as much on its sampling strategy as on its depth or width. Since the underlying spaces are continuous, new training points can be drawn at each epoch, often at negligible extra computational cost. 
This is especially appealing on compact manifolds or on dense coordinate patches where random or quasi-random sampling is easy. 
The network then acts as a smooth interpolant over the entire domain, not merely as a table of values on a mesh. 
For high-order geometric quantities this mesh-free feature is particularly useful, as once the architecture is trained, one can evaluate the learned metric, curvature, or immersion and any autodifferentiated derivatives at arbitrary new points. 
The three applications below exploit this common philosophy in distinct ways: first through a two-patch atlas for Einstein metrics, then through a global conformal-factor model, and finally through a neural immersion whose training directly minimises a curvature energy.

%%%%%%%%%%%%%%%%%%%%%%%%%%%%%%%%%%%%%%%%%%%%%%%%%%%%
\section{Geometric Applications}\label{sec:geo}
The three case studies below highlight complementary ways in which PINNs can be adapted to geometry. The Einstein-metric problem is naturally atlas-based and emphasises patch consistency. The Nirenberg problem is naturally expressed by a global scalar unknown on a fixed manifold and emphasises interpretable post-processing. The Willmore problem treats the immersion itself as the learnable object and shifts the emphasis from solving a PDE residual to minimising a geometric energy directly. 
Taken together they illustrate a useful lesson: the main design freedom is not only the neural network itself, but which geometric quantity one chooses to represent and which constraints are enforced by architecture rather than by loss.

\subsection{Einstein Metrics on Spheres} % use ../ainstein/
The AInstein framework of \cite{Hirst:2025seh} studies Einstein metrics on spheres by combining a two-patch Atlas architecture with a purely geometric training objective. The target equation is
\begin{equation}
Ric(g)=\lambda g, \qquad \lambda\in\{+1,0,-1\},
\label{eq:einstein_eq}
\end{equation}
and the manifolds considered are $S^n$ for $n=2,3,4,5$. Rather than working directly with the standard stereographic coordinates on $\mathbb{R}^n$, the construction composes stereographic projection with a radial map to the unit ball. This yields two finite ball patches related by
\begin{equation}
\tilde{x}=\tau(x)=\frac{|x|-1}{|x|(|x|+1)}x,
\label{eq:einstein_transition}
\end{equation}
which makes both sampling and visualisation substantially cleaner. On each patch a neural subnetwork predicts a lower-triangular vielbein $L$, from which the metric is reconstructed as $g=L^TL$; this guarantees a symmetric metric and its positive definiteness by construction.

The differential content of the problem is then enforced through the standard coordinate expressions
\begin{equation}
\Gamma_{ij}^{k}=\frac{1}{2}g^{k\ell}(\partial_i g_{j\ell}+\partial_j g_{i\ell}-\partial_\ell g_{ij}),
\qquad
R_{jk}=\partial_i\Gamma_{jk}^{i}-\partial_j\Gamma_{ki}^{i}+\Gamma_{ip}^{i}\Gamma_{jk}^{p}-\Gamma_{jp}^{i}\Gamma_{ik}^{p},
\label{eq:einstein_ricci}
\end{equation}
evaluated by automatic differentiation of the network output. 
Training uses a weighted sum of Einstein residuals on the two patches, then overlap losses impose $g^{(1)}=J^Tg^{(2)}J$ on the overlap region, and finiteness losses penalise numerically degenerate outputs. 
The loss weights used in \cite{Hirst:2025seh} are $(f_1,f_2,f_3)=(1,10,1)$, reflecting the importance of patch compatibility in constructing a global geometry.

\begin{figure}[t]
\centering
\begin{subfigure}[t]{0.48\textwidth}
\centering
\includegraphics[width=\textwidth]{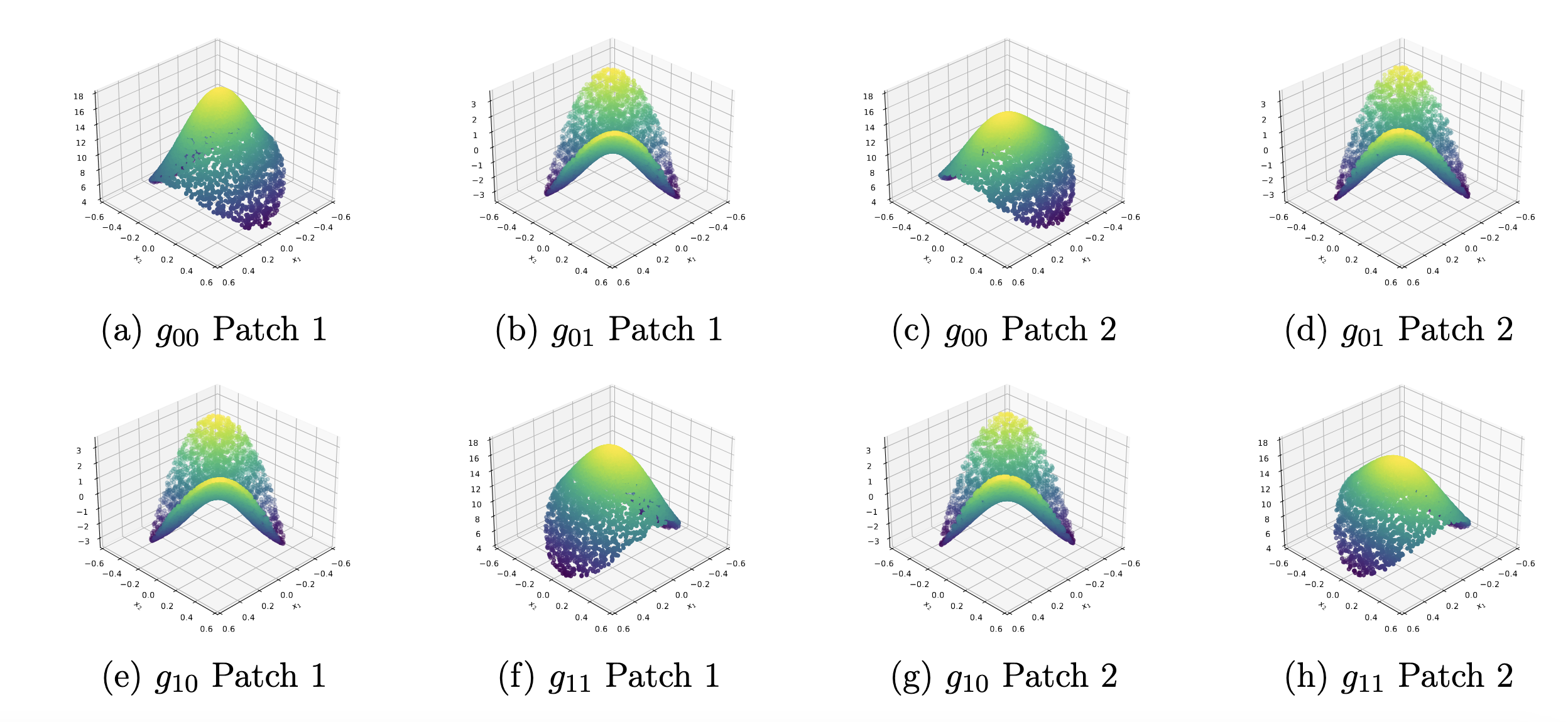}
\caption{Learnt metric components on the $S^2$ patches \\ for $\lambda=+1$.}
\label{fig:einstein_metric_fig}
\end{subfigure}
\hfill
\begin{subfigure}[t]{0.48\textwidth}
\centering
\includegraphics[width=\textwidth]{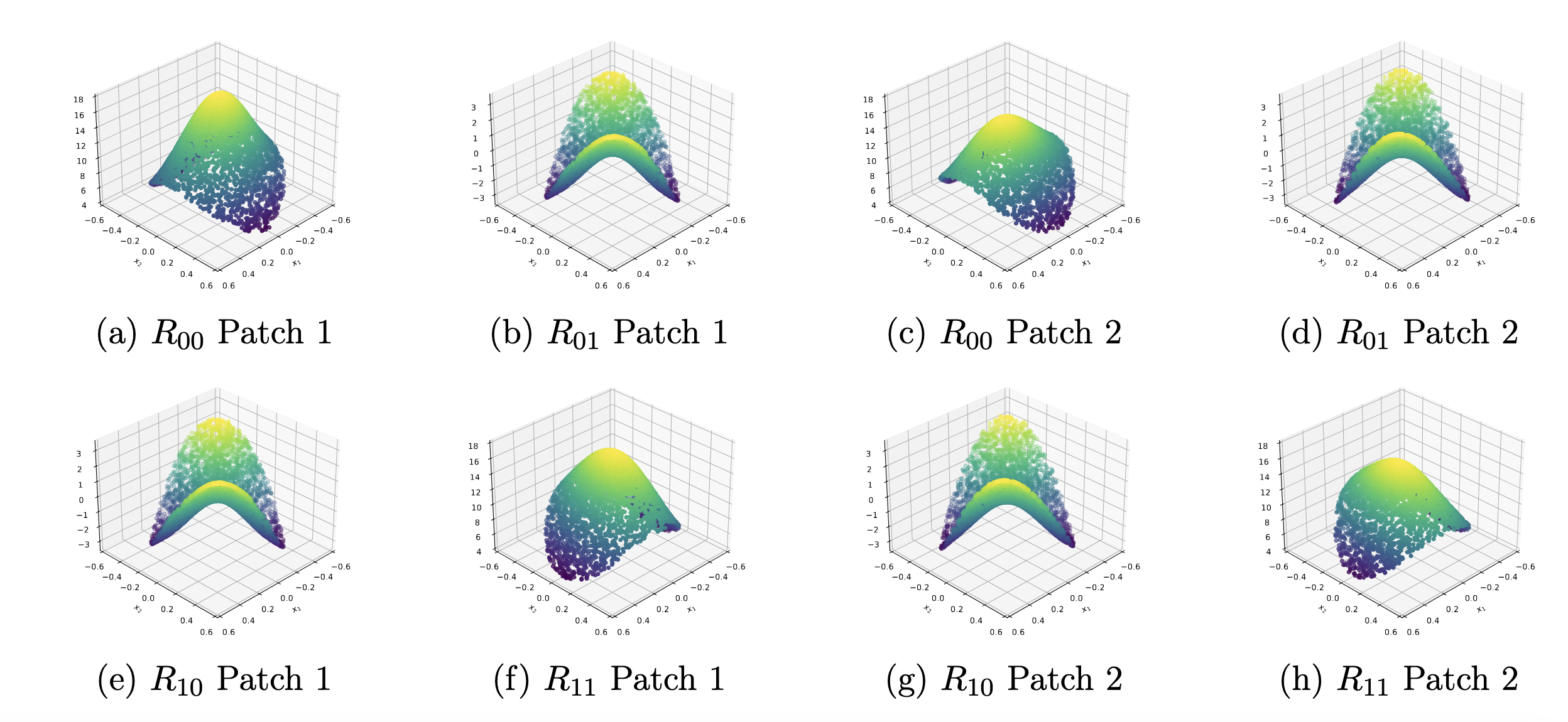}
\caption{Equivalent Ricci tensor components on the  the $S^2$ patches for $\lambda=+1$.}
\label{fig:einstein_ricci_fig}
\end{subfigure}
\caption{Representative two-dimensional positive-curvature output from the AInstein atlas model. 
The $(x,y)$ directions are the $2d$ coordinates of the respective patch, and the $z$ direction is the value of that tensor component. The close agreement between metric and Ricci components corroborates the strong learning of the $Ric(g)=g$ condition directly.}
\label{fig:einstein_results}
\end{figure}

The resulting numerical evidence, with test loss values in Table \ref{tab:einstein_losses}, is intriguing. For the known round metrics with $\lambda=+1$, the semi-supervised PINN attains global test losses comparable to, and in several dimensions slightly better than, a supervised baseline trained directly on the analytic round metric; corroborated by visualisations given in Figure \ref{fig:einstein_results}. For $\lambda=\{0,-1\}$ in dimensions $2$ and $3$, where compact spheres cannot admit such Einstein metrics \cite{Besse:1987pua}, the losses remain large. 
The open and exploratory cases are dimensions $4$ and $5$, where existence of Einstein metrics with $\lambda=\{0,-1\}$ is unknown, and in this work the equivalent high-loss behaviour persists, providing novel numerical evidence against Ricci-flat and hyperbolic Einstein metrics on $S^4$ and $S^5$.

\begin{table}[t]
\centering
\begin{tabular}{|c|ccc!{\vrule width 1.2pt}c|}
\hline
Dimension & $\lambda=+1$ & $\lambda=0$ & $\lambda=-1$ & Supervised $\lambda=+1$ \\
\hline
2 & $0.083\pm0.023$ & $2.881\pm0.113$ & $4.364\pm0.093$ & $0.096\pm0.013$ \\
3 & $0.151\pm0.027$ & $5.560\pm0.160$ & $8.641\pm0.183$ & $0.195\pm0.020$ \\
4 & $0.150\pm0.018$ & $8.494\pm0.121$ & $14.928\pm1.317$ & $0.248\pm0.024$ \\
5 & $0.244\pm0.039$ & $10.810\pm0.185$ & $18.798\pm2.024$ & $0.518\pm0.063$ \\
\hline
\end{tabular}
\caption{Global test losses reported in \cite{Hirst:2025seh} for semi-supervised Einstein-metric learning on $S^n$, compared with a supervised round-metric baseline. The $\lambda=+1$ column tracks the known round solution well, while the $\lambda=0,-1$ columns remain far larger.}
\label{tab:einstein_losses}
\end{table}

Beyond the investigations into the open problems of $\lambda=\{0,-1\}$ in dimensions $\{4,5\}$, the main methodological contribution is that the metric itself is learned with minimal prior symmetry assumptions. 
The atlas formalism is therefore not tied to spheres: the same strategy extends to more complicated manifolds, provided one can choose workable local coordinates and transition maps. 
In that sense AInstein is a prototype for geometric PINNs on spaces represented with atlas methods traditional and central to differential geometric constructions, and compatible where global coordinates are either inconvenient or unavailable.

\subsection{Scalar Curvatures on $S^2$} % use ../nirenberg/
\begin{figure}[t]
\centering
\begin{subfigure}[t]{0.48\textwidth}
\centering
\includegraphics[width=\textwidth]{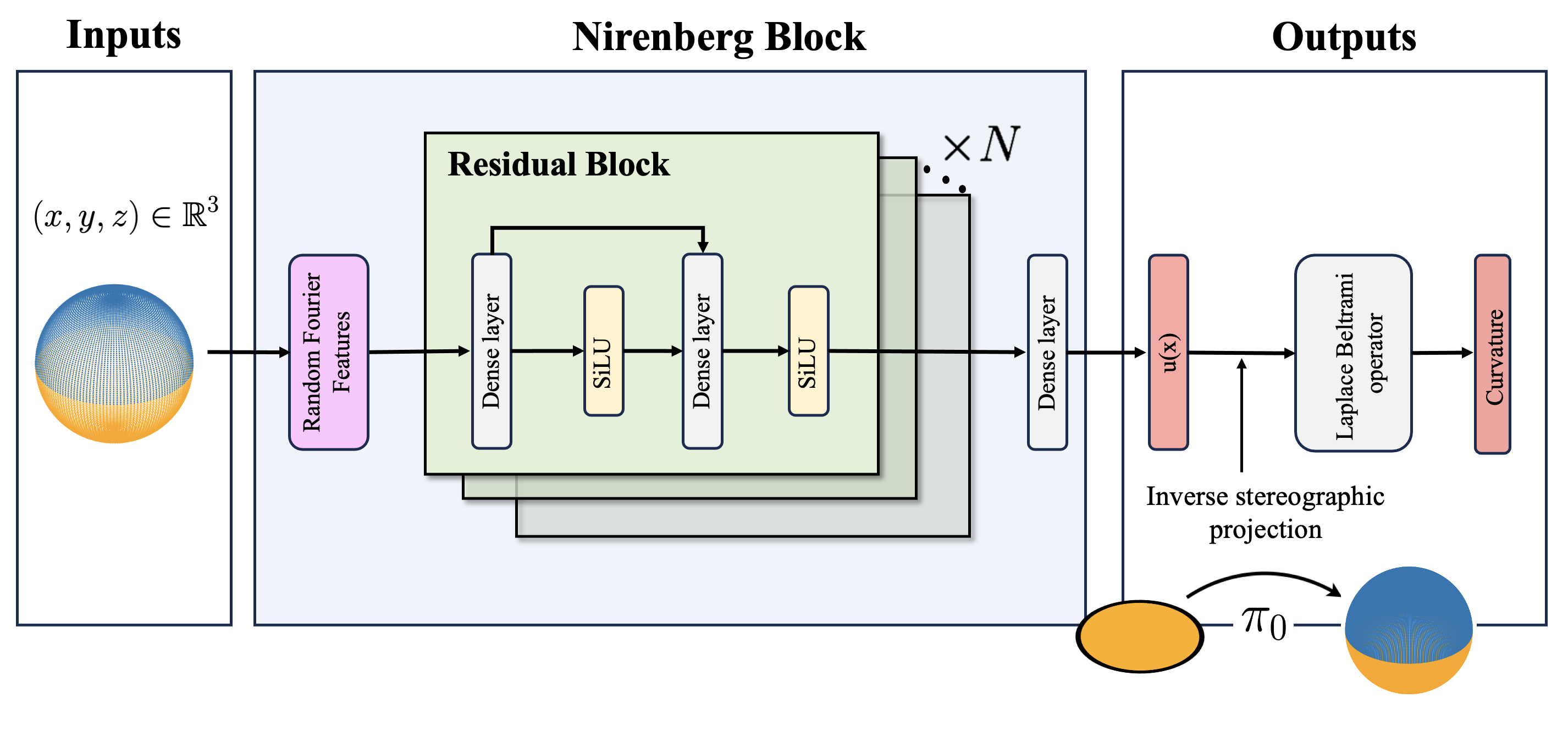}
\caption{Global architecture for learning the conformal factor on $S^2$.}
\label{fig:nirenberg_architecture}
\end{subfigure}
\hfill
\begin{subfigure}[t]{0.48\textwidth}
\centering
\includegraphics[width=\textwidth]{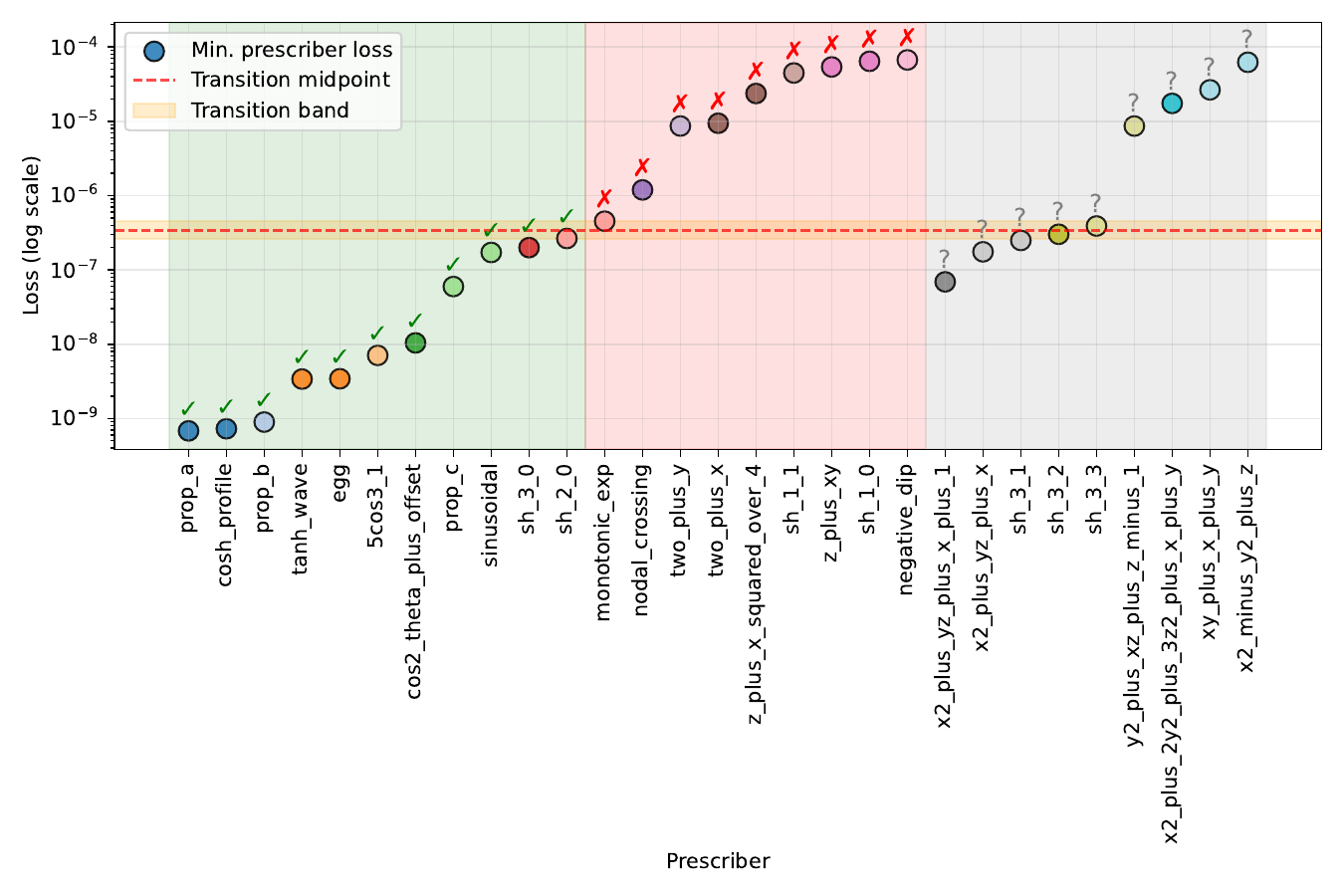}
\caption{Minimum losses across all tested prescribers, separating known realisable, obstructed, and unresolved cases.}
\label{fig:nirenberg_losses}
\end{subfigure}
\caption{This PINN architecture combines a global conformal-factor ansatz with a global optimisation diagnostic. The loss plot shows the same empirical separation between known realisable and obstructed cases that is later used to assess unresolved prescribers.}
\label{fig:nirenberg_overview}
\end{figure}

The work \cite{Cortes:2026kfx} extends this PINN approach to address the Nirenberg problem, originally posed by Nirenberg \cite{Nirenberg1953}, and stated as:\\
\textit{For a given a smooth function $K:S^2\to\mathbb{R}$, does there exist a metric conformal to the round metric $g_0$ whose Gaussian curvature is $K$?}\\
Writing $g=e^{2u}g_0$, since all metrics in $2d$ are conformally equivalent, the problem becomes the semilinear elliptic equation
\begin{equation}
1-\Delta_{g_0}u = K e^{2u}\,,
\label{eq:nirenberg_pde}
\end{equation}
and the problem reduces to finding a suitable $u$.
This problem is especially well suited to a global Embedding architecture, whose design is shown in Figure \ref{fig:nirenberg_architecture}, embedding the $S^2$ sphere in $\mathbb{R}^3$.
Additionally, since the unknown $u$ is a scalar conformal factor, rather than a tensor field, it does not transform between patches or hence require patchwise gluing. 
The paper uses a neural network to represent $u$ directly on $S^2$, with derivatives computed through automatic differentiation and the residual of \eqref{eq:nirenberg_pde} used as the core training signal.

The analytic structure of the Nirenberg problem plays an important role in model design and validation. 
A solution must satisfy the Gauss--Bonnet constraint and the Kazdan--Warner identity \cite{KW1974curvature,KW1975scalar}, while variational existence theory is governed by Moser's functional
\begin{equation}
J(u)=\int_{S^2}\bigl(|\nabla u|^2+2u\bigr)\,d\nu_0-\log\int_{S^2}Ke^{2u}\,d\nu_0,
\label{eq:nirenberg_moser}
\end{equation}
and by the minimax methods developed further in \cite{ChangYang1987}. 
In the computational setting these identities become more than background theory, they are diagnostic tools. 
Low PDE residuals accompanied by small Gauss--Bonnet violations provide evidence that the learned $u$ is not merely fitting local data, but is capturing a globally coherent conformal metric.

The numerical results recover the expected distinction between realisable and obstructed prescribed curvatures. 
For benchmark cases with known realisability, the network reaches losses in the $10^{-7}$ to $10^{-10}$ range, whereas known non-realisable examples remain orders of magnitude larger; shown in Figure \ref{fig:nirenberg_losses}. 
This separation is already useful as an exploratory device, but this work also uses this split as a diagnostic tool for functions with unknown realisability, allowing conjectures to be made about realisability for some spherical harmonics and more general prescribers. 
In particular, the experiments suggest that the previously unresolved harmonics $Y_{3,1}$, $Y_{3,2}$, and $Y_{3,3}$ are realisable, as indicated by low loss values, constituting concrete novel conjectures.
The resulting PINN architecture is therefore not just a solver for known PDE instances; it becomes a computational probe for an open existence problem.
Excitingly, later work in \cite{platt_nirenberg}, inspired by the work summarised here, used rigorous computer-assisted methods to prove the realisability of $Y_{3,2}$ as a scalar curvature on $S^2$, directly validating this AI-assisted discovery process in differential geometry. 

To make the learned solutions interpretable, the authors fit a truncated spherical-harmonic expansion
\begin{equation}
\tilde{u}=\sum_{\ell\leq L}\sum_{m=-\ell}^{\ell}\frac{c_{\ell,m}}{\ell(\ell+1)}Y_{\ell,m}
\label{eq:nirenberg_harmonic_ansatz}
\end{equation}
to the network output. This produces explicit approximate conformal factors and hence explicit approximate metrics. For the prescribers whose conformal factors are already known from spectral-pair constructions, the learned coefficients reproduce the expected sparse structure accurately; for previously unresolved cases, they provide the first interpretable candidates for the underlying geometry.

From the perspective of PINN methodology, the Nirenberg study shows how classical analytic obstructions and modern computational optimisation can complement one another. The geometric analysis narrows the search space and supplies non-trivial validation criteria, while the neural model provides a flexible global ansatz that can be interrogated after training. This is precisely the sort of interaction between numerical experimentation and analytic conjecture that makes PINNs attractive in differential geometry.

\begin{figure}[t]
\centering
\begin{subfigure}[t]{0.32\textwidth}
\centering
\includegraphics[width=\textwidth]{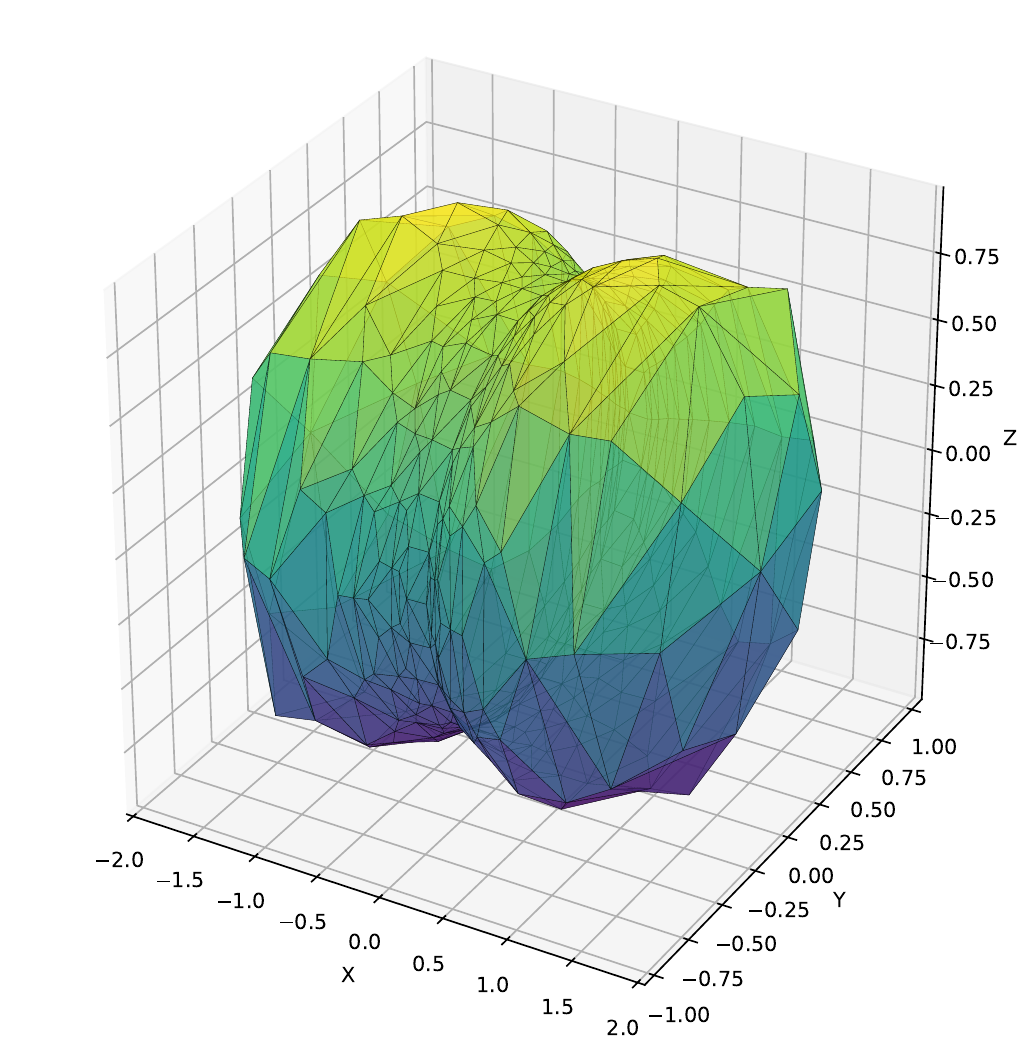}
\caption{$Y_{2,0}$ (realisable)}
\label{fig:nirenberg_embed_sh20}
\end{subfigure}
\hfill
\begin{subfigure}[t]{0.32\textwidth}
\centering
\includegraphics[width=\textwidth]{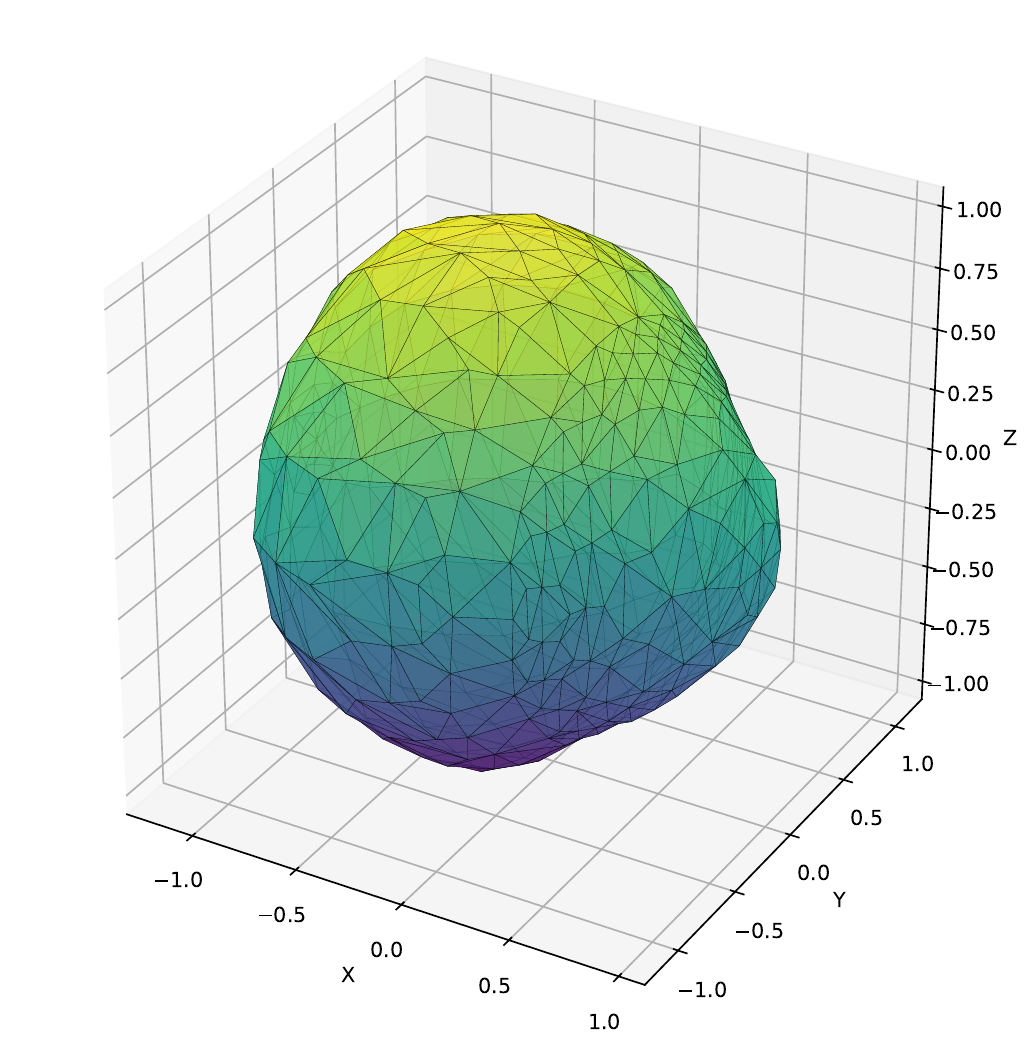}
\caption{$Y_{3,0}$ (realisable)}
\label{fig:nirenberg_embed_sh30}
\end{subfigure}
\hfill
\begin{subfigure}[t]{0.32\textwidth}
\centering
\includegraphics[width=\textwidth]{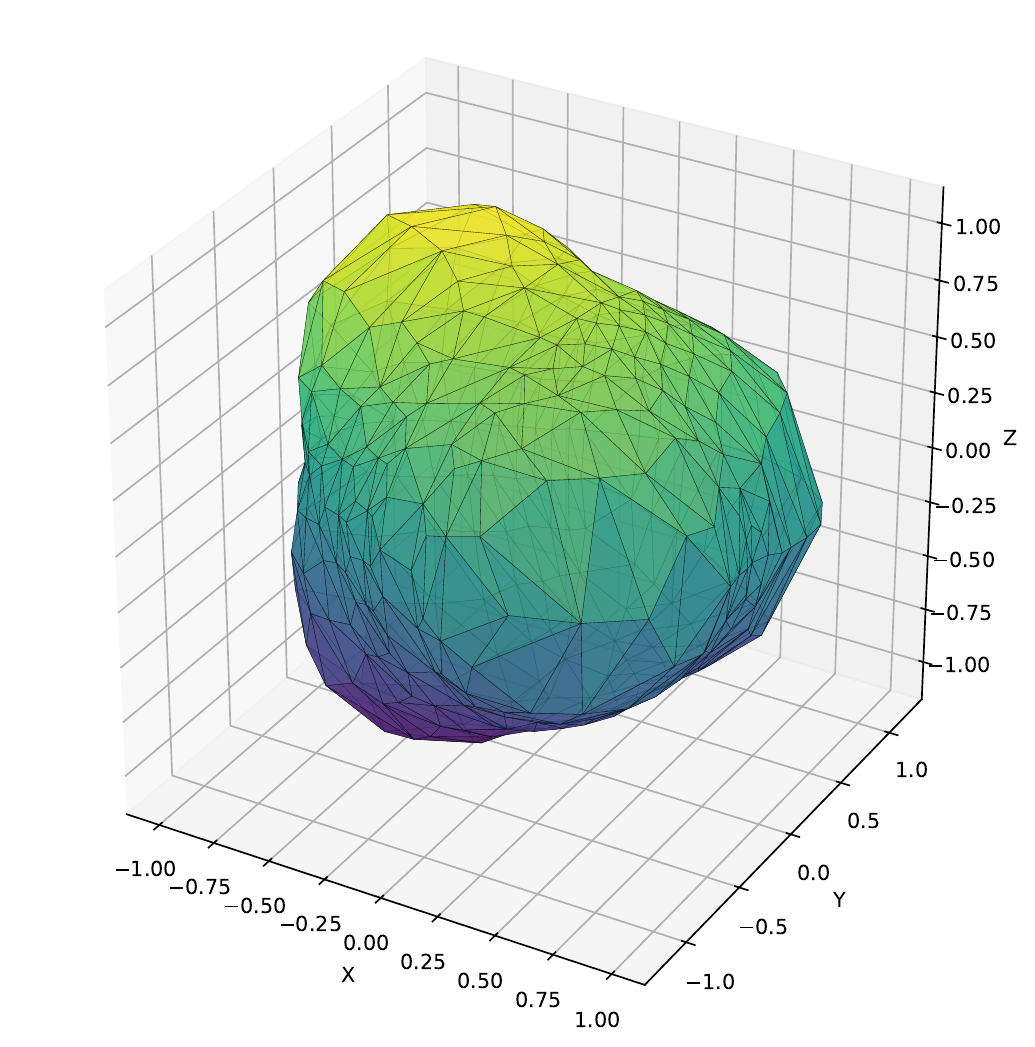}
\captionsetup{justification=centering}
\caption{$Y_{3,2}$ (unknown, \\ conjectured realisable)}
\label{fig:nirenberg_embed_sh32}
\end{subfigure}
\caption{MDS embeddings of learned metrics for three spherical-harmonic prescribers. The first two are benchmark cases with known realisability, while the third illustrates the geometry of the case $Y_{3,2}$ with conjectured realisability.}
\label{fig:nirenberg_metric_visualisations}
\end{figure}

The learned metrics can also be visualised directly through MDS embeddings of the geodesic distance data induced by the recovered conformal factors. In Figure~\ref{fig:nirenberg_metric_visualisations}, the benchmark harmonics $Y_{2,0}$ and $Y_{3,0}$ display the smooth deformations expected of these known realisable cases, while the $Y_{3,2}$ metric sits in the regime that was originally unresolved which loss results in \cite{Cortes:2026kfx} then motivated conjectured realisability (later proved in \cite{platt_nirenberg}). This geometric visualisation therefore complements the loss-based diagnostics by giving an interpretable picture of the candidate metrics themselves.

\subsection{Minimal Willmore Surfaces in $\mathbb{R}^3$} % use ../willmore/
The third application, developed in \cite{Hirst:2026qwi}, shifts from learning a metric or conformal factor to learning the immersion itself. Given a smooth immersion $\varphi:\Sigma\to\mathbb{R}^3$ from a fixed $2d$ topological surface, the first and second fundamental forms may be written in the $2d$ local coordinates $(u,v)$ of the fundamental domain as
\begin{equation}
\mathrm{I}=E\,du^2+2F\,du\,dv+G\,dv^2,
\qquad
\mathrm{II}=L\,du^2+2M\,du\,dv+N\,dv^2,
\end{equation}
with coefficients for $\mathrm{I}$ as
\begin{equation}
E=\langle \varphi_u,\varphi_u\rangle,
\qquad
F=\langle \varphi_u,\varphi_v\rangle,
\qquad
G=\langle \varphi_v,\varphi_v\rangle,
\end{equation}
where $\varphi_i$ represents partial derivatives with respect to coordinate $i$. Then writing
\begin{equation}
\hat{n}=\frac{\varphi_u\times\varphi_v}{|\varphi_u\times\varphi_v|},
\end{equation}
\begin{equation}
L=\langle \varphi_{uu},\hat{n}\rangle,
\qquad
M=\langle \varphi_{uv},\hat{n}\rangle,
\qquad
N=\langle \varphi_{vv},\hat{n}\rangle,
\end{equation}
for the $\mathrm{II}$ coefficients.
The first fundamental form is the induced metric. Its area element is $dA=\sqrt{EG-F^2}\,du\,dv$. The second fundamental form records how the immersion bends in $\mathbb{R}^3$.
From these the Willmore energy introduced by Willmore \cite{Willmore1965} is defined in terms of the mean curvature, $H$, as
\begin{equation}
\mathcal{W}(\varphi)=\int_\Sigma H^2\,dA,
\qquad
H=\frac{EN-2FM+GL}{2(EG-F^2)}.
\label{eq:willmore_energy}
\end{equation}

Instead of discretising the surface by a mesh and then approximating curvature numerically, the paper represents $\varphi$ by a smooth neural map from a fundamental domain. Automatic differentiation then supplies the derivatives needed to compute $H$ and the area element exactly at sampled points.

The architecture is built on genus-dependent fundamental domains and uses spectral input maps to enforce the required identifications before optimisation. For genus $0$, on $[0,2\pi]\times[0,\pi]$, the inputs are low-order real spherical harmonics, which automatically become independent of the azimuthal angle at the poles and therefore close the sphere smoothly. For genus $1$, and for each torus chart used in genus $2$, the doubly periodic domain $[0,2\pi]^2$ is passed through Fourier features
\begin{equation}
\Phi(u,v)=\bigl(\sin ku,\cos ku,\sin kv,\cos kv\bigr)_{k=1}^{N},
\end{equation}
so points differing by integer multiples of $2\pi$ in either direction have identical network inputs and the torus periodicity is enforced exactly. In the genus-$2$ case these periodic torus charts are punctured by removing a small disc, where the surrounding annuli are then glued by additional regularity losses to form the target higher-genus surface.

This leads naturally to a Monte Carlo PINN loss,
\begin{equation}
\mathcal{L}_{\mathcal{W}}(\varphi_\theta)=|\Omega|\,\frac{1}{B}\sum_{i=1}^{B}H_i^2\sqrt{E_iG_i-F_i^2},
\label{eq:willmore_mc}
\end{equation}
augmented by regularity losses preventing metric degeneration and, in the genus-$2$ case, these gluing losses that enforce $C^2$ compatibility between two punctured-torus charts. The method may be viewed as a neural Willmore flow: gradient descent in parameter space continuously deforms the immersion so as to reduce the curvature energy while preserving the intended topology through the choice of parameter domain. 
General initialisation of the architecture gives a completely random (often regularly self-intersecting) immersion.
Therefore to maintain the correct genus, supervised pretraining was used to seed the model with deliberately non-optimal reference surfaces with the correct topology which could then be evolved. These initial surfaces were an ellipsoid for genus $0$, a simple torus for genus $1$, and two displaced punctured tori for genus $2$.

\begin{figure}[t]
\centering
\begin{subfigure}[t]{0.32\textwidth}
\centering
\includegraphics[width=\textwidth]{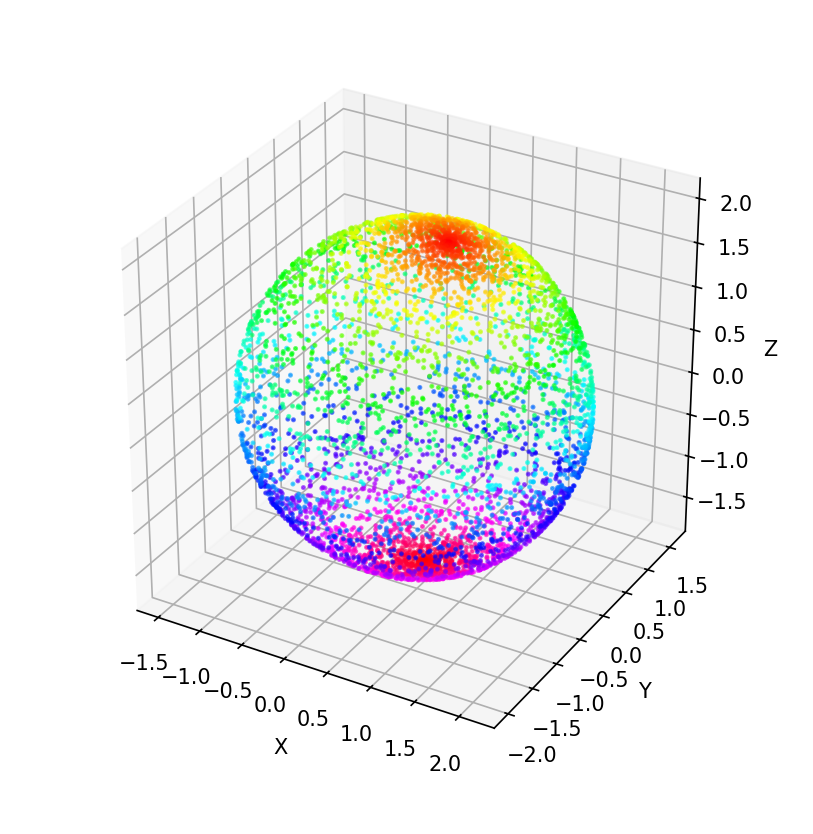}
\caption{Genus $0$ final surface.}
\label{fig:willmore_g0}
\end{subfigure}
\hfill
\begin{subfigure}[t]{0.32\textwidth}
\centering
\includegraphics[width=\textwidth]{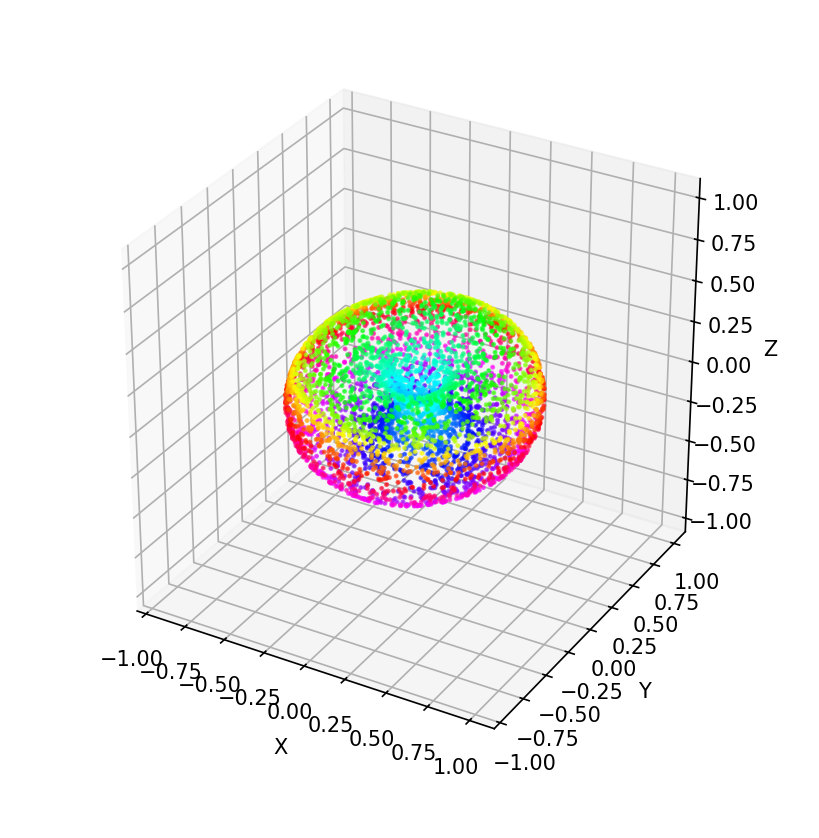}
\caption{Genus $1$ final surface.}
\label{fig:willmore_g1}
\end{subfigure}
\hfill
\begin{subfigure}[t]{0.32\textwidth}
\centering
\includegraphics[width=\textwidth]{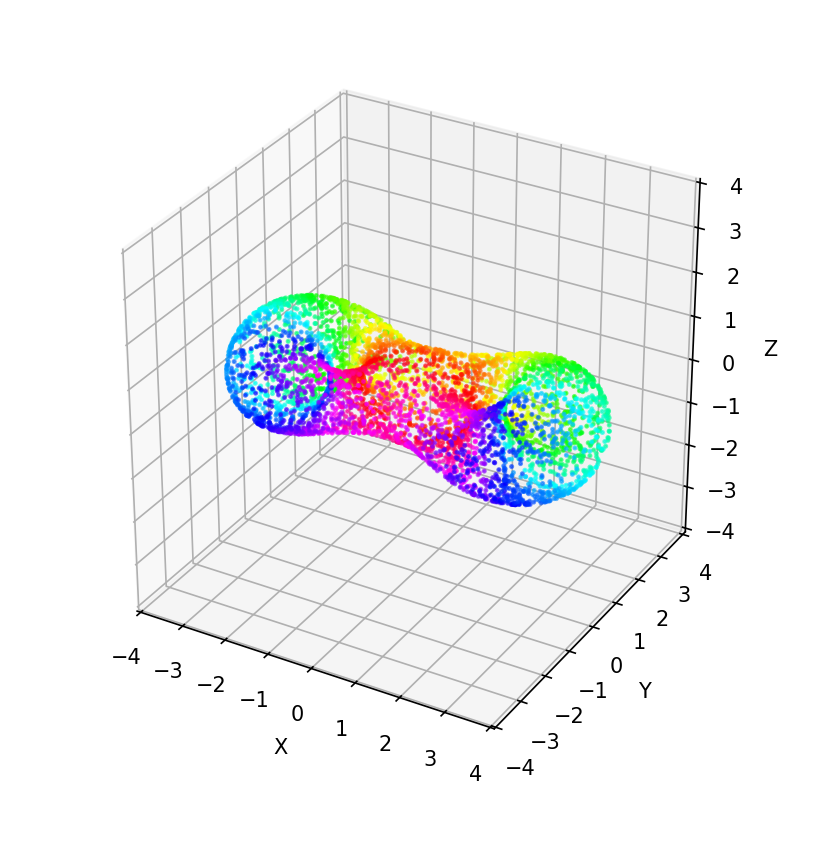}
\caption{Genus $2$ final surface.}
\label{fig:willmore_g2}
\end{subfigure}
\\[0.8em]
\begin{subfigure}[t]{0.32\textwidth}
\centering
\includegraphics[width=\textwidth]{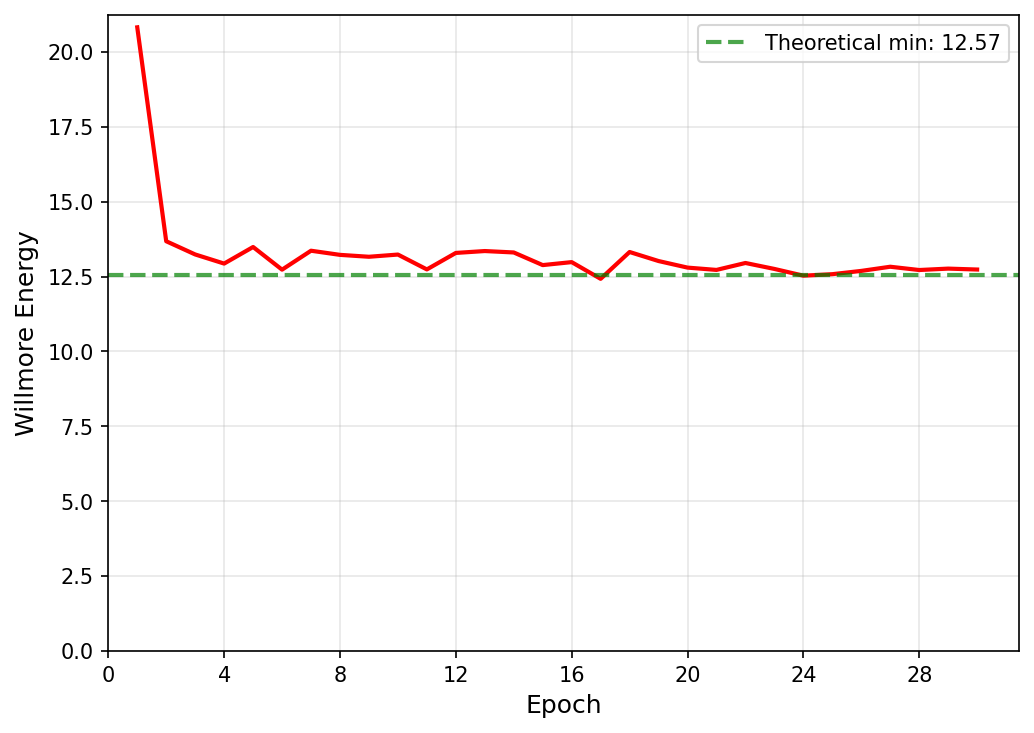}
\caption{Genus $0$ Willmore loss.}
\label{fig:willmore_loss_g0}
\end{subfigure}
\hfill
\begin{subfigure}[t]{0.32\textwidth}
\centering
\includegraphics[width=\textwidth]{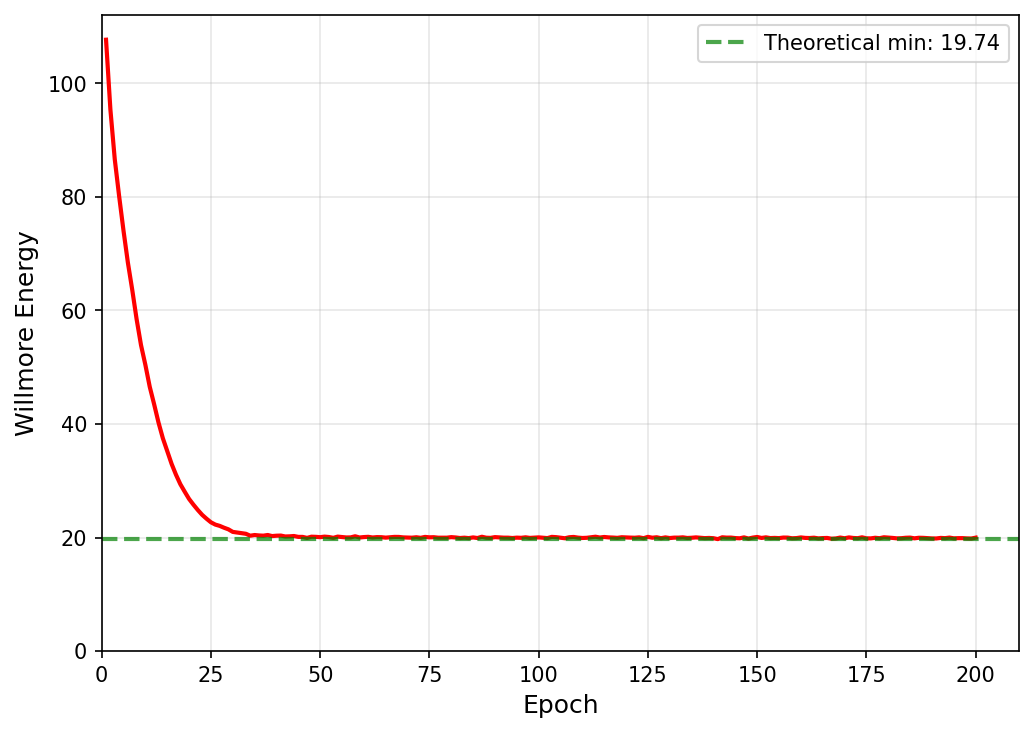}
\caption{Genus $1$ Willmore loss.}
\label{fig:willmore_loss_g1}
\end{subfigure}
\hfill
\begin{subfigure}[t]{0.32\textwidth}
\centering
\includegraphics[width=\textwidth]{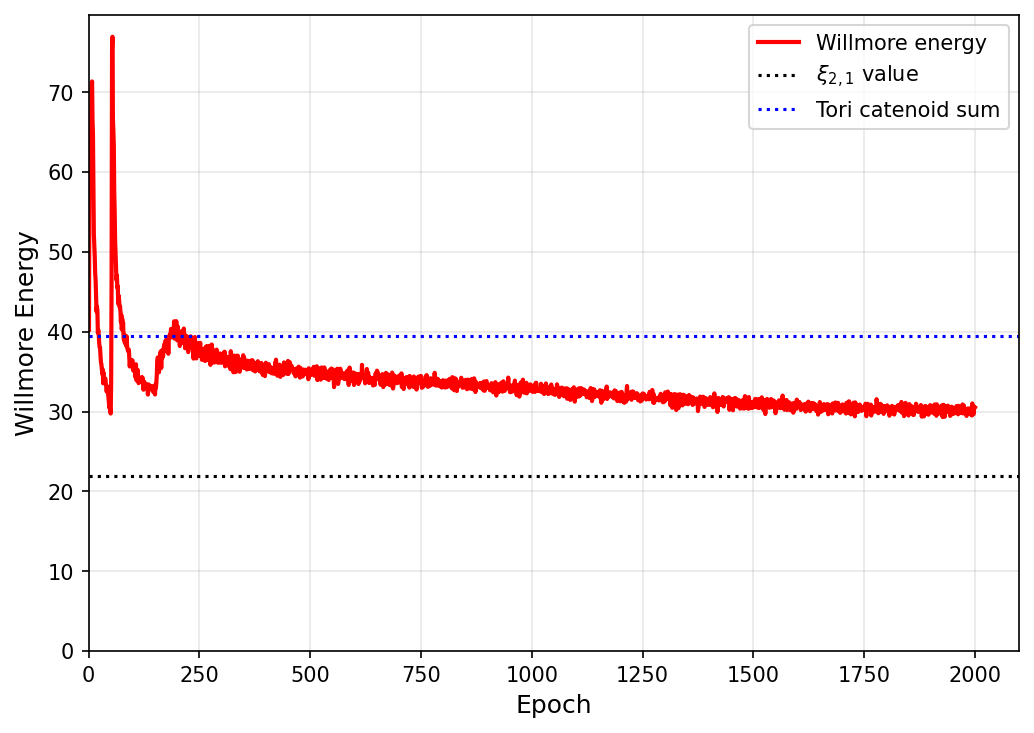}
\caption{Genus $2$ Willmore loss.}
\label{fig:willmore_energy_trace}
\end{subfigure}
\caption{Neural Willmore flow outputs from \cite{Hirst:2026qwi}. The first row shows the final learned embeddings for genuses $0$, $1$, and $2$; the second row shows the corresponding Willmore-losses over training.}
\label{fig:willmore_results}
\end{figure}

\begin{figure}[t]
\centering
\begin{subfigure}[t]{0.19\textwidth}
\centering
\includegraphics[width=\textwidth]{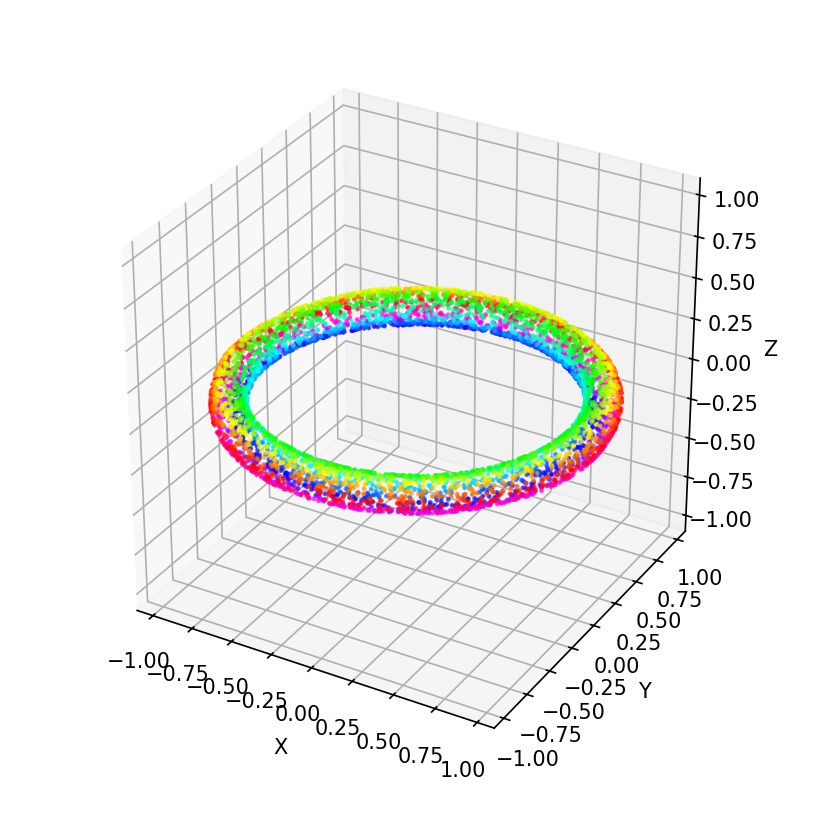}
\captionsetup{justification=centering}
\caption{Epoch $0$ \\ \hspace*{3em} $\widehat{\mathcal{W}}=118.97$}
\label{fig:willmore_g1_epoch0}
\end{subfigure}
\hfill
\begin{subfigure}[t]{0.19\textwidth}
\centering
\includegraphics[width=\textwidth]{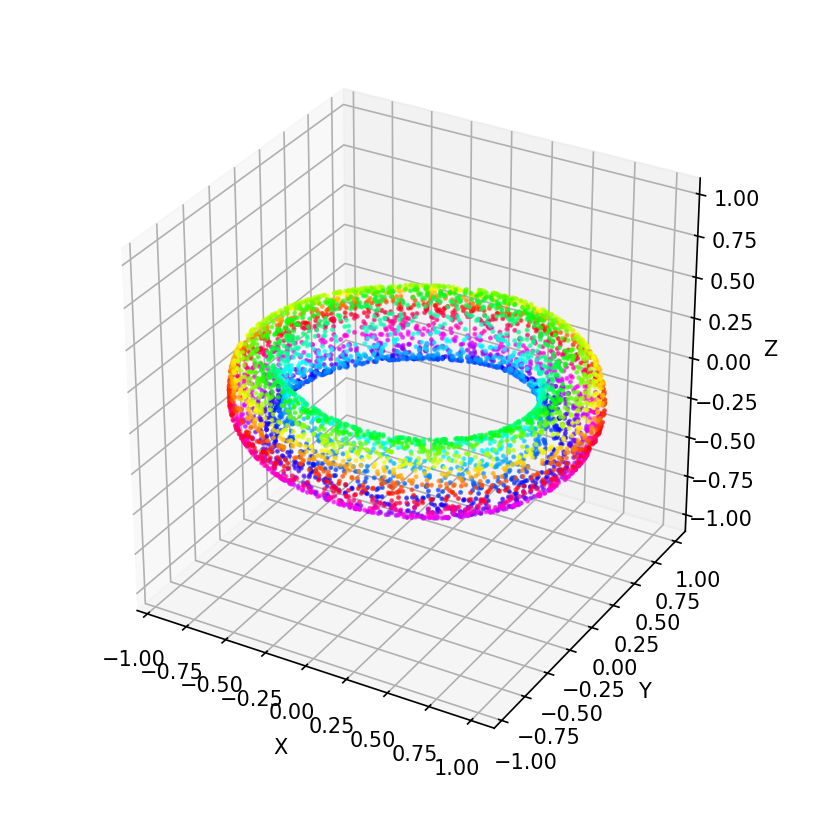}
\captionsetup{justification=centering}
\caption{Epoch $10$ \\ \hspace*{2em} $\widehat{\mathcal{W}}=50.42$}
\label{fig:willmore_g1_epoch10}
\end{subfigure}
\hfill
\begin{subfigure}[t]{0.19\textwidth}
\centering
\includegraphics[width=\textwidth]{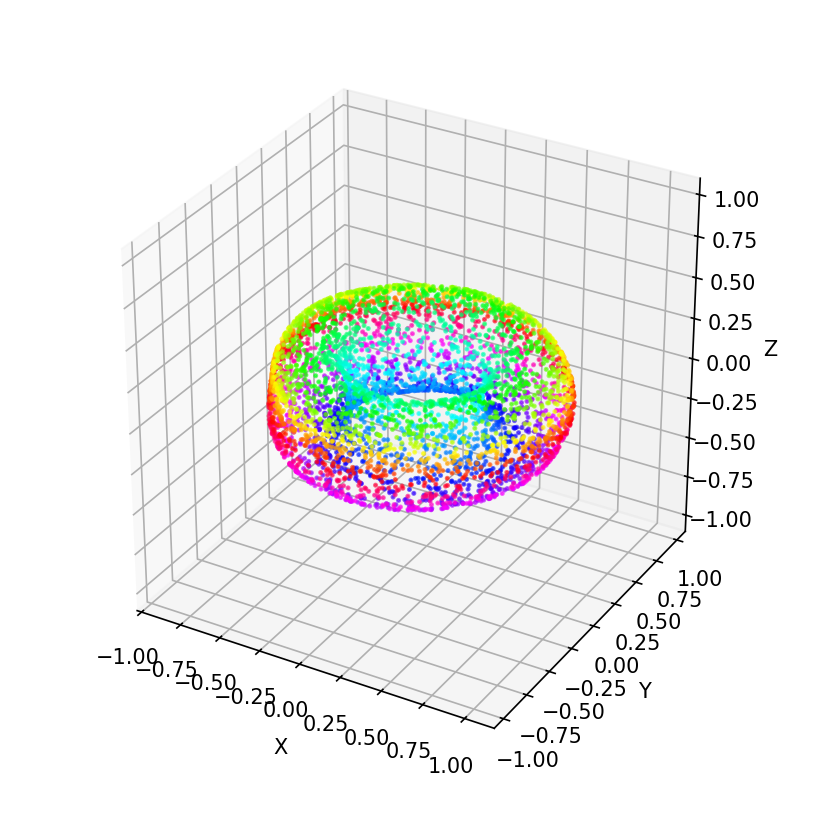}
\captionsetup{justification=centering}
\caption{Epoch $20$ \\ \hspace*{2em} $\widehat{\mathcal{W}}=26.80$}
\label{fig:willmore_g1_epoch20}
\end{subfigure}
\hfill
\begin{subfigure}[t]{0.19\textwidth}
\centering
\includegraphics[width=\textwidth]{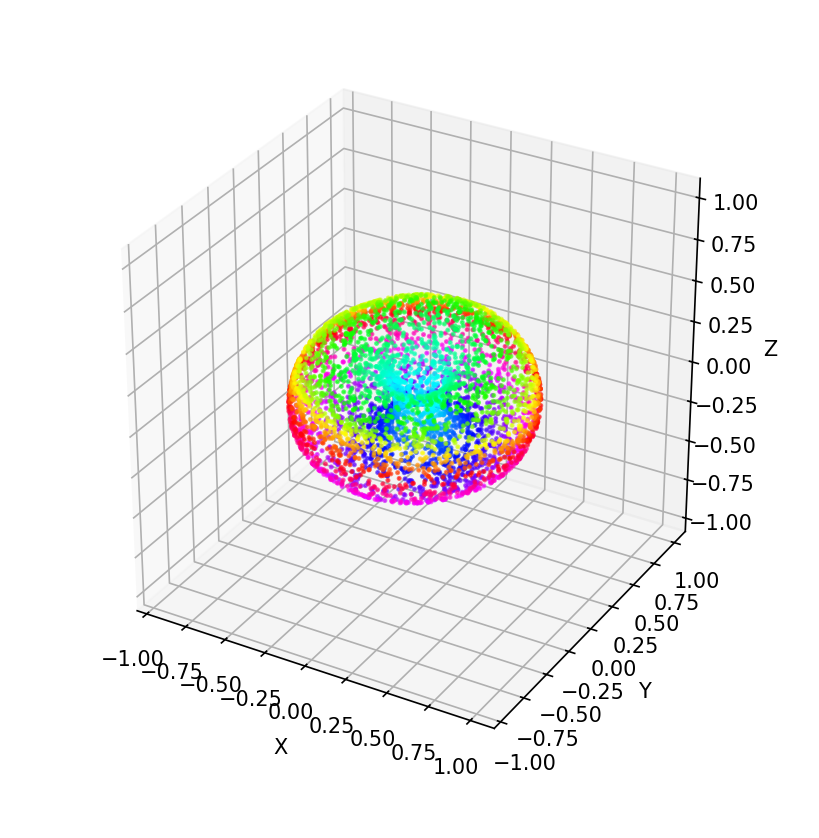}
\captionsetup{justification=centering}
\caption{Epoch $100$ \\ \hspace*{1.5em} $\widehat{\mathcal{W}}=20.03$}
\label{fig:willmore_g1_epoch100}
\end{subfigure}
\hfill
\begin{subfigure}[t]{0.19\textwidth}
\centering
\includegraphics[width=\textwidth]{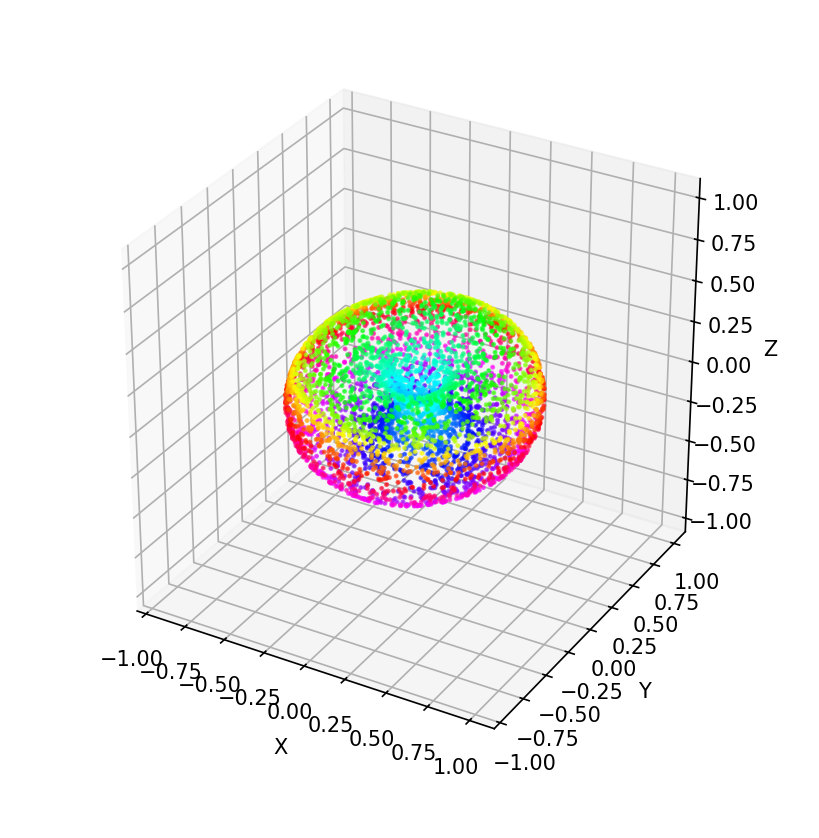}
\captionsetup{justification=centering}
\caption{Epoch $200$ \\ \hspace*{1.5em} $\widehat{\mathcal{W}}=19.98$}
\label{fig:willmore_g1_epoch200}
\end{subfigure}
\caption{Five snapshots from the genus-$1$ learning process, taken from the embedding sequence reported in \cite{Hirst:2026qwi}. The optimisation rapidly lowers the Willmore energy, converging to the Clifford-torus minimum.}
\label{fig:willmore_g1_evolution}
\end{figure}

Final learnt surfaces, as well as the respective Willmore loss curves, are given in Figure \ref{fig:willmore_results}. The validation cases are genus $0$ and genus $1$, where the analytic minima are known. For genus $0$ the network evolves an initial ellipsoid toward the round sphere and reaches a final Monte Carlo estimate $\widehat{\mathcal{W}}\approx 12.73$, close to the exact value $4\pi\approx 12.57$. For genus $1$ it evolves a torus toward the Clifford torus and reaches $\widehat{\mathcal{W}}\approx 19.98$, close to the sharp minimum $2\pi^2\approx 19.74$ proved by Marques and Neves \cite{MarquesNeves2014}. These two cases are important because the architecture has no explicit bias toward those canonical surfaces beyond the topology and regularity encoded in the setup.
A more detailed view of the genus-$1$ optimisation is given in Figure~\ref{fig:willmore_g1_evolution}, showing the surface evoltion over training.

The genus-$2$ experiment is the most exploratory part of the paper. Guided by Kusner's conjecture that the minimiser should be the stereographic projection of the Lawson surface $\xi_{2,1}$ \cite{Kusner1989,Lawson1970}, the network is initialised from two punctured tori and trained with staged gluing losses. The resulting immersed surface achieves $\widehat{\mathcal{W}}\approx 30.19$, significantly below the naive catenoid-bridge heuristic value $4\pi^2\approx 39.48$, though still above the Lawson benchmark. This does not resolve the genus-$2$ problem, but it demonstrates that a PINN-style immersion ansatz can navigate a genuinely non-trivial search space and produce quantitatively meaningful candidates in an open geometric problem.

Conceptually, this example shows that PINNs need not be tied to residual formulations of PDEs. When a geometric problem is fundamentally variational, the loss can be the energy itself. The network then serves as a smooth parameterisation of a class of admissible objects, and optimisation becomes a new numerical flow on that class, useful throughout geometry.

%%%%%%%%%%%%%%%%%%%%%%%%%%%%%%%%%%%%%%%%%%%%%%%%%%%%
\section{Conclusions}\label{sec:conc}
PINN architectures are a particularly natural fit for differential geometry because the objects one wants to learn are smooth, the conditions defining them are differential, and the underlying spaces are continuous such that they support dense resampling during training. The main design choice is not simply which optimiser or hidden-layer width to use, but which geometric quantity is represented by the network and which constraints are encoded architecturally rather than left to the loss. Across the examples reviewed here, that principle leads to three distinct but compatible strategies: atlas-based metric learning for Einstein geometry, global scalar-field learning for prescribed curvature on $S^2$, and direct immersion learning for Willmore minimisation.

The three case studies also show different ways in which neural methods can contribute to geometry. In the Einstein problem they reproduce known positive-curvature solutions and provide new negative evidence in an open setting. In the Nirenberg problem they separate solvable from obstructed prescribers and suggest new candidates in unresolved harmonic families. In the Willmore problem they recover the round sphere and Clifford torus, then move into open territory at genus $2$ with quantitatively non-trivial low-energy surfaces. These are not replacements for analytic arguments, but they are already useful as conjectural support devices, as sources of explicit candidates, and as a way to organise geometric constraints into a trainable computational framework.

The broader lesson is that PINNs become more compelling, not less, as the geometry becomes richer. More complicated geometries offer more interesting coordinate systems, compatibility conditions, conserved quantities, and variational structures to exploit. The recent progress summarised here suggests that neural architectures, when designed with those structures in mind, can become a practical part of the modern computational toolkit for differential geometry.

%%%%%%%%%%%%%%%%%%%%%%%%%%%%%%%%%%%%%%%%%%%%%%%%%%%%
\section*{Acknowledgements}
The author wishes to thank \{Tancredi Schettini Gherardini, Alexander G. Stapleton\}, \{Gianfranco Cortés, Maria Esteban-Casadevall, Yueqing Feng, Jonas Henkel, Tancredi Schettini Gherardini, Alexander G. Stapleton\}, and \{Henrique N. Sá Earp, Tomás S.R. Silva\} for collaboration on the respective works summarised here; and acknowledges support from São Paulo Research Foundation (FAPESP) grant 2024/18994-7. \\
This paper summaries work presented at the workshop "Recent Progress in Computational String Geometry" 26w5653 in January 2026, hosted by Banff International Research Station at the Chennai Mathematical Institute; the author thanks the workshop organisers for the invitation to attend, speak, and contribute this proceedings.

%%%%%%%%%%%%%%%%%%%%%%%%%%%%%%%%%%%%%%%%%%%%%%%%%%%%
\newpage
\linespread{0.9}\selectfont
%\addcontentsline{toc}{section}{References}
\bibliographystyle{utphys}
\bibliography{references}

\end{document}